\documentclass[
12pt,
a4paper]{article}
\newcommand { \Proof} {\noindent {\bf Proof: }}
 \newtheorem {theorem} {Theorem} 
 \newtheorem {prop} {Proposition} 
 \newtheorem {lemma} [prop] {Lemma}
 \newtheorem {definition} [prop] {Definition}
  
\newcommand { \qed} { $  \hskip.2cm \hfill \blacksquare$ \hskip.5cm}

 \newcommand{\R}{{\sf R\hspace*{-0.9ex}\rule{0.15ex}{1.4ex}\hspace*{0.9ex}}}
 \advance \textheight by \topskip
\input epsf
\usepackage{amssymb}
\usepackage[pdftex]{graphicx,color} 
\title { Generic Oval Billiards} 
\author 
{M\'ario Jorge DIAS CARNEIRO, Sylvie OLIFFSON KAMPHORST\\
  and S\^onia PINTO-DE-CARVALHO \\
Departamento de Matem\'atica, ICEx, UFMG \\
30.123-970, Belo Horizonte, Brasil} 
\date {Versao final submetida, 7 de maio 2007} 

\begin {document}
\maketitle 
\begin{abstract}
In this paper we show that, under certain generic conditions, billiards on ovals have only a finite number of periodic orbits, for each period $N$, all non-degenerate and  at least one of them is hyperbolic.
Moreover, the invariant curves of two hyperbolic points are transversal.
We explore these properties to give some dynamical consequences  specially about the dynamics in the instability regions.
\end{abstract}

\section{Introduction}

Let $\alpha$ be an oval, i.e., a planar, closed, regular, simple, oriented counterclockwise, $C^l$ curve, $l\geq 2$, parameterized by $t\in S^1$ and such that its radius of curvature $R(t)>0$.
The billiard problem on $\alpha$ consists in the free motion of a point
particle  on the planar region enclosed by $\alpha$, being reflected
elastically at the impacts on the boundary. The trajectories will be polygonals
on this planar region.

Since the motion is free inside the region, it is completely determined by the points of impacts at $\alpha$, and the direction of motion immediately after each reflection,
defined by the angle $\theta$ between it and the oriented tangent to the boundary at the reflection point.
Therefor, to each oval $\alpha$ is associated a billiard map $T_\alpha$ from the cylinder $S^1\times(0,\pi)$  into itself which to each initial condition $(t_0,\theta_0)$ associates the next impact and direction: $T_\alpha(t_0,\theta_0)=(t_1,\theta_1)$.

The map $T_\alpha$ has some very well known properties (see for instance \cite{bir27}, \cite{kat},  \cite{koz}, \cite{rag}, \cite{tab}): It is a $C^{l-1}$-diffeomorphism, preserving the measure  $d\mu = \frac{ds}{dt}\sin\theta \, d\theta \, dt$, where $s=s(t)$ is the arclength parameter for $\alpha$. It is reversible with respect to the reversing symmetry $H(t,\theta)=(t,\pi-\theta)$, which satisfies $H\circ H=Id$ and $H\circ T_\alpha=T_\alpha^{-1}\circ H$. It is a monotone Twist map.
So, the billiard map defines a discrete reversible conservative bidimensional dynamical system in the cylinder $S^1\times(0,\pi)$. This map has no fixed points but Birkhoff \cite{bir27} proved that it has periodic orbits of any period greater or equal to 2. 

In this paper we show that, under certain generic conditions, billiards on ovals have only a finite number of periodic orbits, for each period $N$, all non-degenerate and  at least one of them is hyperbolic.  Moreover, the invariant curves of two hyperbolic points are transversal. 
The generic existence, in the analytic case, of at least one nondegenerate periodic orbit for each period, was proved by Kozlov and Treschev \cite{koz}. Our results generalize theirs in the sense that we work with curves at least $C^2$, prove the finiteness and study the heteroclinic connections.

Once we have a hyperbolic periodic orbit, which unstable manifold project to $S^1$,
we can construct the instability region associated to it from the closure of the invariant manifold. 
The dynamical object we obtain  is the same as the one constructed by Le Calvez in \cite{lec87} and thus shares the same dynamical properties.

\section{Generic Periodic Orbits}

Billiard maps are a special kind of diffeomorphisms and $C^1$ perturbations of the map may produce diffeomorphisms which are not billiards. 
In order to assure that the perturbed diffeomorphism remains a billiard, we have to perturb the boundary curve instead of the map itself.

As planar rigid motions and homotheties do not change the geometrical features of the curve, they do not change the dynamical characteristics of the associated billiard map and so it is useful to work on the set of curves modulo this equivalence relation.
Let ${\cal C}$ be the set of equivalence classes of $C^2$ ovals.
Given  a representative $\alpha$ of an equivalence class $[\alpha]\in{\cal C}$, the normal bundle $(\alpha(t),\eta(t))$ is $C^1$.
For $\epsilon>0$, a tubular neighbourhood of $\alpha$ is given by
$N_\epsilon(\alpha)=\{\alpha(t)+\lambda\eta(t),\, 0\le t < 1,\, -\epsilon<\lambda<\epsilon \} $.

\begin{definition}
$[\beta]\in{\cal C}$ is $\epsilon$-close to $[\alpha]\in{\cal C}$ if there exist representatives $\alpha$ in $[\alpha]$ and $\beta$ in $[\beta]$ such that the image of $\beta$ is in $N_\epsilon(\alpha)$, and its canonical projection onto the image of $\alpha$ is a diffeomorphism.
\end{definition}
 As a consequence, $\beta$ can be written as 
$\beta(t)=\alpha(t)+\lambda(t)\eta(t)$, with
$\lambda$ at least $C^2$ and periodic.

\begin{definition}
$[\beta]\in{\cal C}$ is $\epsilon$-$C^2$-close to $[\alpha]\in{\cal C}$ if
$[\beta]$ is $\epsilon$-close to $[\alpha]$ and $ ||\lambda||_2   <\epsilon$.
\end{definition}

\begin{prop}\label{prop:baire}
${\cal C}$ is a Baire space.
\end{prop}
\Proof $C^2(S^1,\R^2)$ with the $C^2$-topology is a Baire space. By Sard's Theorem, the subset $I$ of immersions is open and dense in $C^2(S^1,\R^2)$, and so also a Baire space.
For $\alpha\in I$ let $R(t)$ be its radius of curvature on each $t$ and let $index$ be the total curvature divided by $2\pi$. Let $A=\{\alpha \hbox{ such that } R\geq 0 \hbox{ and } index =1\}$ and $B=\{\alpha$ such that $R> 0$ and $index\ =1\}$. 
Note that $B$ is exactly the set of closed, regular, simple, oriented, $C^2$-curves
with positive radius of curvature. 
It is clear that $A$ is closed and  $B$ is open and dense in $A$, so $B$, the set of ovals, is Baire with the $C^2$-topology.
Cutting by the equivalence relation and taking the induced topology, we have that ${\cal C}$ is a Baire space. \qed

As we have showed in \cite{etds}, by perturbing  the classes of curves we get nearby diffeomorphisms in the $C^1$-topology:

\begin{prop}\label{prop:difeo}
If $[\beta]\in{\cal C}$ is $\epsilon$-$C^2$-close to $[\alpha]\in{\cal C}$ then there exist representatives $\alpha$ in $[\alpha]$ and 
$\beta$ in $[\beta]$ such that the associated billiard maps $T_{\alpha}$ and $T_{\beta}$ are $C^1$-close.
\end{prop}

To lighten the notation, in what follows we will denote both the curve $\alpha$ and the equivalence class $[\alpha]$ by $\alpha$, unless where confusion may happen.

\subsection{Periodic Orbits}

Let $\alpha\in {\cal C}$. Since the radius of curvature $R(t)$ is strictly positive, $\alpha$ can be reparameterized by $\varphi\in [0,2\pi)$, the positive oriented angle between the tangent vector $\alpha'(t)$ and a fixed direction (say, the x-axis).
Let $T_\alpha:[0,2\pi)\times (0,\pi)\to [0,2\pi)\times (0,\pi) $ be the associated billiard map.

A point $(\overline \varphi_0,\overline\theta_0)\in [0,2\pi)\times(0,\pi)$ is $n$-periodic if $n$ is the smallest positive integer such that 
$T_\alpha^n(\overline \varphi_0,\overline \theta_0)=(\overline \varphi_0,\overline\theta_0)$. Such a point defines an $n$-periodic orbit 
$$\{(\overline \varphi_0,\overline\theta_0), T_\alpha(\overline \varphi_0,\overline\theta_0)=(\overline \varphi_1,\overline\theta_1),...,T_\alpha^{n-1}(\overline \varphi_0,\overline\theta_0)=(\overline \varphi_{n-1},\overline\theta_{n-1})\}.$$

An $n$-periodic orbit is nondegenerate if each of its points 
$(\overline \varphi_j,\overline\theta_j)$ is a nondegenerate fixed point of $T_\alpha^n$. Otherwise it is degenerate.

Let $U_N$ the set of ovals $\alpha \in \cal C$ such that for each $n\neq 1$
divisor of $N$, the associated billiard map $T_\alpha$ has only nondegenerate $n$-periodic orbits. 

\begin{prop} If $\alpha\in U_N$ then for each $n\neq 1$
divisor of $N$, $T_\alpha$ has only a finite number of $n$-periodic orbits.
\end{prop}

\Proof 
Let  $\{(\overline \varphi_0,\overline\theta_0), (\overline \varphi_1,\overline\theta_1),...,(\overline \varphi_{n-1},\overline\theta_{n-1})\}$ be an $n$-periodic orbit of $T_{\alpha}$.
Then $\alpha(\overline \varphi_0)$, $\alpha(\overline \varphi_1)$,..., $\alpha(\overline \varphi_{n-1})$ are the vertices of a polygon inscribed in $\alpha$.
Let $\beta_i$ be its internal angles and $\gamma_i=\overline\theta_i$ if $\overline\theta_i\leq\pi/2$ or $\gamma_i=\pi-\overline\theta_i$ if $\overline\theta_i\geq\pi/2$.
As $\sum_{i=0}^{n-1}\beta_i\leq (n-2)\pi$, $2\gamma_i+\beta_i=\pi$ and then
$\sum_{i=0}^{n-1}\gamma_i\geq\pi$. 
So there is at least one $j$ such that 
$\pi/n\leq \overline\theta_j\leq (n-1)\pi/n$. 

This means that if $\alpha\in U_N$ then each $n$-periodic orbit of $T_\alpha$, with $n\neq 1$ divisor of $N$, has at least one point in the compact cylinder $[0,2\pi)\times [\pi/N,(N-1)\pi/N]$. 

Since all the points of $T^N_\alpha$ are nondegenerate, there is only a finite number of them on the compact cylinder and so, only a finite number of $n$-periodic orbits.\qed

For diffeomorphisms $C^1$ on $[0,2\pi)\times(0,\pi)$, with the $C^1$-topology, having only a finite number of nondegenerate fixed points on a compact subset is an open property. 
Taking the restriction to billiard maps, we get 

\begin{prop} $U_N$ is open on $ \cal C$. \label{prop:open} \end{prop}

By Birkhoff's Theorem, any $\alpha\in \cal C$ has periodic orbits of any period. So $\alpha$ will be outside $U_N$ if it has (finitely or infinitely many) degenerate $n$-periodic orbits for $n\neq 1$, divisor of $N$. 
The following lemma provides the basic tool for proving the density of $U_N$.

\begin{lemma}\label{lem:perturb}
Suppose that ${\cal O}=\{(\overline \varphi_0,\overline\theta_0),
(\overline \varphi_1,\overline\theta_1),...,
(\overline \varphi_{n-1},\overline\theta_{n-1})\}$ is a degenerate $n$-periodic orbit for the billiard map $T_\alpha$ associated to a curve $\alpha\in \cal C$. Then there is a curve $\tilde\alpha\in \cal C$, $C^2$-close to $\alpha$, such that $\cal O$ is a nondegenerate $n$-periodic orbit for the associated billiard map 
$T_{\tilde\alpha}$.
\end{lemma}

\Proof 
Let $x_j=R(\varphi_j)\sin\theta_j$ and 
$l_j=||\alpha(\varphi_j)-\alpha(\varphi_{j+1})||$. The
Jacobian matrix of $T_\alpha$ at $(\varphi_j,\theta_j)$ is given by:
$$DT_\alpha (\varphi_j, \theta_j)=
\frac{1}{x_{j+1}}
 \left( 
\begin{array}{cc} 
l_j-x_j & l_j \\
 l_j- x_j- x_{j+1} &
 l_j- x_{j+1}
\end{array}\right ).$$
and, for any point $(\overline \varphi_j,\overline \theta_j)$ of  ${\cal O}$, det$\left(DT^n_\alpha(\overline \varphi_j,\overline \theta_j)\right)=1$. It follows that ${\cal O}$ is degenerate if $\mbox{tr}\left(DT^n_\alpha(\overline \varphi_j,\overline \theta_j)\right)=\pm 2$, for any point of the orbit.

We have that
$$\mbox{tr}\left(DT^n_\alpha(\overline \varphi_0,\overline\theta_0)\right)=
\mbox{tr}\left(DT_\alpha{(\overline \varphi_{n-1},\overline \theta_{n-1})}...
DT_\alpha{(\overline \varphi_0,\overline \theta_0)}\right)=
\mbox{tr}\left(A_{n-1}A_{n-2}...A_{1}A_{0}\right),$$
where 
$$A_{j}=
\left( 
\begin{array}{cc} 
\displaystyle \frac{l_j}{x_j}-1 & 
\displaystyle \frac{l_j}{x_j} \\
\displaystyle \frac{l_j}{x_j}-\frac{x_{j+1}}{x_j}-1  &
\displaystyle \frac{l_j}{x_j}- \frac{x_{j+1}}{x_j}
\end{array}\right )
$$
Let us isolate the terms that depend on, for instance, $x_1$.
Only $A_{1}$ and $A_{0}$ have entries with $x_1$ and 
\begin{eqnarray}
A_1 \, A_0 &=&
\frac{2}{x_1} \,
\left( 
\begin{array}{cc} 
 \frac{l_1 \, l_0}{x_0}-l_1 & 
 \frac{l_1 \, l_0}{x_0} \\
 \frac{l_1 \,l_0}{x_0}- \frac{l_0 \, x_{2}}{x_0}- l_1 - + x_2  &
 \frac{l_1 \, l_0}{x_0}- \frac{l_0 \, x_{2}}{x_0}
\end{array}\right )
- \left( 
\begin{array}{cc} 
 \frac{l_1+l_0}{x_0}-1 & 
 \frac{l_1+l_0}{x_0} \\
 \frac{l_1+l_0}{x_0}- \frac{x_{2}}{x_0}-1  &
 \frac{l_1+l_0}{x_0}- \frac{x_{2}}{x_0}
\end{array}\right ) \nonumber \\
&=& \frac{1}{x_1}B_1+ C_{1,0} \label{eqn:c}  
\end{eqnarray}
Then 
\begin{eqnarray*}
\mbox{tr}\left(DT^n_\alpha(\overline \varphi_0,\overline\theta_0)\right)&=&
\mbox{tr}\left[A_{n-1}A_{n-2}...A_2\left(\frac{1}{x_1}B_1+C_{1,0}\right)\right]\\
&=&\frac{1}{x_1}\mbox{tr}\left(A_{n-1}A_{n-2}...A_2B_1\right)+
\mbox{tr}\left(A_{n-1}A_{n-2}...A_2C_{1,0}\right)\\
&=&\frac{b_1}{x_1}+c_1
\end{eqnarray*}
where neither $b_1$ nor $c_1$ depend on $x_1$.

If $b_1\neq 0$, let $I$ be an interval such that $\overline \varphi_1\in I$ and $\overline \varphi_j\notin I$ for $j\neq 1$. 
Let $\tilde\alpha(\varphi)=\alpha (\varphi)+\lambda (\varphi)\eta(\varphi)$, where $\lambda$ is a $C^2$ periodic function satisfying 
$\lambda(\varphi)=0$ if $\varphi\notin I$, $\lambda(\overline \varphi_1)=\lambda'(\overline \varphi_1)=0$,
$\lambda''(\overline \varphi_1)\neq 0$ and  $||\lambda ||_2$  small enough to guarantee that $\tilde R$, the radius of curvature of $\tilde\alpha$, is strictly positive. 

The perturbed curve $\tilde\alpha$ and the original one $\alpha$ coincide, except on a neighbourhood of $\alpha(\overline\varphi_1)$ and, at $\overline\varphi_1$, they have a contact of order one.
So the polygonal trajectory that corresponds, on the billiard
table, to the periodic orbit is unchanged and $\cal O$ is also a $n$-periodic orbit for $T_{\tilde\alpha}$.  

Moreover $\tilde l_j=||\tilde\alpha(\overline\varphi_j)-\tilde\alpha(\overline\varphi_{j+1})||=l_j, \forall j$, the radius of curvature  
$\tilde R(\overline\varphi_j)=R(\overline\varphi_j)$ if $j\neq 1$ and 
$\tilde R(\overline\varphi_1)=R(\overline\varphi_1)\left(1-\frac{\lambda''(\overline\varphi_1)}{R(\overline\varphi_1)}\right)^{-1}$.

Then $$\mbox{tr}\left(DT^n_{\tilde\alpha}(\overline \varphi_0,\overline\theta_0)\right)= \frac{b_1}{\tilde x_1 }+c_1=
\frac{b_1}{\tilde R(\overline\varphi_1)\sin\overline\theta_1 }+c_1=
\frac{b_1}{x_1}(1-\frac{\lambda''(\overline\varphi_1)}{R(\overline\varphi_1)})+c_1$$ 
and we can choose $\lambda$ as small as we want such that 
$\mbox{tr}\left(DT^n_{\tilde\alpha}(\overline \varphi_0,\overline\theta_0)\right)\neq \pm 2$ and $\cal O$ is a nondegenerate $n$-periodic orbit for $T_{\tilde\alpha}$.

If $b_1=0$ then 
$$\mbox{tr}\left(DT^n_\alpha(\overline \varphi_0,\overline\theta_0)\right)
= \mbox{tr} \left(A_{n-1}A_{n-2}...A_2C_{1,0}\right) $$
with
\begin{eqnarray*}
A_2 \, C_{1,0}= \frac{1}{x_2}B_2+C_{2,0}
\end{eqnarray*}
where $C_{2,0}$ has the same form of $A_0$, replacing 
$x_1$ by $x_3$ and $l_0$ by $l_2+l_1+l_0$. 

Then $$\mbox{tr}\left(DT^n_\alpha(\overline \varphi_0,\overline\theta_0)\right)
=\frac{b_2}{x_2}+c_2$$
where neither $b_2$ nor $c_2$ depend on $x_1$ and $x_2$.
If $b_2\neq 0$ then we can make the normal perturbation on a neighbourhood of 
$\overline\varphi_2$ as above. 

If $b_2=0$ we continue the process until finding a $b_i \ne 0$ and then making the normal perturbation at $\overline\varphi_i$ or to end up with all $b_i$'s $=0$, for $i=1..n-1$ in which case, as $x_0=x_n$, we will have
\begin{eqnarray*}
\mbox{tr}\left(DT^n_\alpha(\overline \varphi_0,\overline\theta_0)\right)
&=&
\mbox{tr}\left(C_{n-1,0}\right )\\
&=&
(-1)^{n-1} \, \mbox{tr} 
 \left( 
\begin{array}{cc} 
 \frac{l_{n-1}+\ldots+l_1+l_0}{x_0}-1 & 
 \frac{l_{n-1}+\ldots+l_1+l_0}{x_0} \\
 \frac{l_{n-1}+\ldots+l_1+l_0}{x_0}- 2  &
 \frac{l_{n-1}+\ldots+l_1+l_0}{x_0}- 1
\end{array}\right ) 
\\
&=&  \pm 2 \left (  \frac{l_{n-1}+\ldots+l_1+l_0}{x_0}-1 \right )
\end{eqnarray*}

As $l_{n-1}+\ldots+l_1+l_0$ is the perimeter of the polygonal trajectory and then different from 0, we can perform the normal perturbation on a neighbourhood of $\overline\varphi_0$ as above  concluding the proof of the lemma.\qed

\begin{prop} $U_N$ is dense on $ \cal C$. \label{prop:dense} \end{prop}
\Proof
Given  $\alpha\notin U_N$, let ${\cal P}_\alpha$ be the set of all fixed points of $T_\alpha^N$ and $\Pi:[0,2\pi)\times (0,\pi)\mapsto [0,2\pi)$ be the projection on the first coordinate.

As $\alpha\notin U_N$, there is a $(\overline\varphi_0,\overline\theta_0)\in {\cal P}_\alpha$ such that ${\cal O}(\overline\varphi_0,\overline\theta_0)$ is a degenerate $n$-periodic orbit of $T_\alpha$ for $n\neq 1$, divisor of $N$.
By lemma \ref{lem:perturb} we can find a curve $\alpha_1$, close to $\alpha$, such that ${\cal O}(\overline\varphi_0,\overline\theta_0)$ is a non-degenerate $n$-periodic orbit of $T_{\alpha_1}$ and so there are intervals $J_0,J_1,...,J_{n-1}\subset [0,2\pi)$ such that $\Pi\left( T_{\alpha_1}^j(\overline\varphi_0,\overline\theta_0)\right)$ is the unique point of $\Pi({\cal  P}_{\alpha_1})$ in $J_j$.

As $S^1$ is compact, we can construct after a finite number of steps, a curve $\alpha_2$, as close as we want to $\alpha$, such that $\Pi({\cal  P}_{\alpha_2})=\{\psi_0,\psi_1,...,\psi_k\}$ and  intervals ${\cal J}_0,{\cal J}_1,...,{\cal J}_k \subset [0,2\pi)$, with $\cup {\cal J}_j=[0,2\pi)$, and such that each $\psi_j$ is the unique point of $\Pi({\cal  P}_{\alpha_2})$ in ${\cal J}_j$ and so each fiber $\Pi^{-1}(\psi_j)$ contains at least one nondegenerate fixed point and maybe other degenerate or nondegenerate fixed points.

If $T_{\alpha_2}^N$ has only non-degenerate fixed points, the problem is solved.
Otherwise,
let $(\psi_0,\theta_0)$ and $(\psi_0,\xi_0)$ be two different fixed points of $T_{\alpha_2}^N$, with $(\psi_0,\theta_0)$ non-degenerate and $(\psi_0,\xi_0)\in [0,2\pi)\times [\pi/N,(N-1)\pi/N]$ degenerate. 

If  $(\psi_0,\xi_0)=T_{\alpha_2}^j(\psi_0,\theta_0)$ then they are both non-degenerate, so suppose that they do not belong to the same orbit. Let $n$ and $m$ be the periods of the orbits, respectively. Then $\mbox{tr}\left(DT^n_{\alpha_2}(\psi_0,\theta_0)\right)\neq\pm 2$ and $\mbox{tr}\left(DT^m_{\alpha_2}(\psi_0,\xi_0)\right)=\pm 2$.

As in the proof of lemma~\ref{lem:perturb}, there is an $l$ such that 
$\mbox{tr}\left(DT^m_{\alpha_2}(\psi_0,\xi_0)\right)= \frac{b_l}{x_l} + c_l$ with $b_l \ne 0$ and
then a normal perturbation in the interval containing $\psi_l$ will produce 
 a curve $\beta$, $C^2$-close to $\tilde\alpha$, such that both of $(\psi_0,\theta_0)$
 and $(\psi_0,\xi_0)$ are nondegenerate fixed points of $T_\beta^N$.

To finish the proof of the proposition, we remark that $[\pi/N,(N-1)\pi/N]$ is a closed interval. So, after a finite number of perturbations we can construct a curve in $U_N$, arbitrarily close to the given $\alpha \notin U_N$.
 \qed

Taking the intersection of the open and dense sets $U_N$s and remembering that ${\cal C}$ is a Baire space it follows that

\begin{theorem}
Generically, for billiards on ovals, there is only a finite number of periodic orbits, for each period $N$, and they are all non-degenerate.
\end{theorem}

\subsection{Intersection of stable and unstable manifolds}

Let $\alpha$ be a$C^l$ oval, $T_\alpha$ the associated billiard map and $F_\alpha:\R\times (0,\pi)\mapsto \R\times (0,\pi)$ a lift of $T_\alpha$. If ${\cal O}(\overline \varphi_0,\overline \theta_0)=\{(\overline \varphi_0,\overline \theta_0),(\overline \varphi_1,\overline \theta_1),...,(\overline \varphi_{n-1},\overline \theta_{n-1})\}$ is an $n$-periodic orbit  of $T_\alpha$ then there is an $m$, $1\leq m<n$ such that
$F^n_\alpha(\overline \varphi_0,\overline \theta_0)=(\overline \varphi_0+2m\pi,\overline \theta_0)$. We will say that ${\cal O}(\overline \varphi_0,\overline \theta_0)$ is an $(m,n)$-periodic orbit.

Given $n$ and $m$, $1\leq m<n$, let 
$$G_{m,n}( \psi_0, \psi_1,,..., \psi_{n-1})=-||\alpha(\psi_0)-\alpha(\psi_1)||-||\alpha(\psi_1)-\alpha(\psi_2)||-...-||\alpha(\psi_{n-1})-\alpha(\psi_0+2m\pi)||.$$
${\cal O}(\overline \varphi_0,\overline \theta_0)=\{(\overline \varphi_0,\overline \theta_0),(\overline \varphi_1,\overline \theta_1),...,(\overline \varphi_{n-1},\overline \theta_{n-1})\}$ is a non-degenerate $(m,n)$-periodic orbit if and only if 
$( \overline\varphi_0, \overline\varphi_1$, $\ldots$, $\overline\varphi_{n-1})$ is a non-degenerate singularity of $G_{m,n}$.
By the Mackay-Meiss criterion \cite{mac}, non-degenerate minima of $G_{m,n}$ correspond to hyperbolic orbits and non-degenerate maxima to elliptic or reverse hyperbolic ones.

If $\alpha\in U_N$ then all the $n$-periodic orbits, $n\neq 1$ divisor of $N$, are nondegenerate singularities for the appropriate $G_{m,n}$ and then are hyperbolic or elliptic. Actually, for each $n$ and each $m$ such that gcd$(m,n)=1$, the first step in the proof of Birkhoff's Theorem is the existence of a global minimum of $G_{m,n}$. 
So there is at least one hyperbolic $n$-periodic for $T_\alpha$. The stability of the other periodic orbits is strongly related to the geometry of the oval $\alpha$, and they can even be all hyperbolic, like the 2-periodic orbits in the examples given in  \cite{etds} or \cite{kozl}.

Let $\{(\overline \varphi_0,\overline\theta_0),(\overline \varphi_1,\overline\theta_1),...,(\overline \varphi_{n-1},\overline\theta_{n-1})\}$ be one $n$-periodic hyperbolic orbit. Associated to each point $(\overline \varphi_j,\overline\theta_j)$ there are two $C^{l-1}$ invariant curves 
$$W^u(\overline \varphi_j,\overline\theta_j)=\{(\varphi,\theta)
\hskip.2cm \mbox{such that}\hskip.2cm T_\alpha^{in}(\varphi,\theta)\to
(\overline \varphi_j,\overline\theta_j)\hskip.2cm \mbox{as}\hskip.2cm i\to-\infty\}$$
and
$$W^{s}(\overline \varphi_j,\overline\theta_j)=\{(\varphi,\theta)
\hskip.2cm \mbox{such that}\hskip.2cm T_\alpha^{in}(\varphi,\theta)\to
(\overline \varphi_j,\overline\theta_j)\hskip.2cm \mbox{as}\hskip.2cm i\to\infty\}$$
called, respectively, unstable and stable curves of 
$(\overline \varphi_j,\overline\theta_j)$.

A point $(\varphi_0,\theta_0)$ is hetero(homo)clinic if $(\varphi_0,\theta_0)\in 
W^u(\overline \varphi_j,\overline\theta_j)\cap W^{s}(\overline \varphi_k,\overline\theta_k)$ for $k\neq j$ ($k=j$).

For the hyperbolic periodic orbits corresponding to global minima of $G_{m,n}$, Bangert's results \cite{ban} assures the existence of hetero and homoclinic points. Other hyperbolic periodic orbits may have only homoclinic points, as the 2-periodic orbit plotted in the figure bellow. 


We don't know if, for a generic oval $\alpha$ there exists a hyperbolic periodic orbit without hetero or homoclinic points. In particular, $W^u(\overline \varphi_j,\overline\theta_j)\cap W^{s}(\overline \varphi_k,\overline\theta_k)=\emptyset$ is an open property for $C^1$-diffeomorphisms, and so will be open for billiards on ovals.

An heteroclinic or homoclinic point $(\varphi_0,\theta_0)$ is called transversal (tangent) if the invariant stable and unstable curves meet transversally (tangentially) at $(\varphi_0,\theta_0)$. Transversal intersection of stable and unstable curves is also an open property for $C^1$-diffeomorphisms, and so it will be open for billiards on ovals. In the next lemma we will prove that any billiard with a tangent hetero(homo)clinic point can be approached by billiards with a transversal one. We use the techniques introduced in \cite{lev} and used in \cite{don}.

\begin{lemma}\label{lem:levallois}
Let $\alpha\in U_N$ and $\{(\overline \varphi_0,\overline \theta_0),...,(\overline \varphi_{n-1},\overline \theta_{n-1})\}$ be an $n$-periodic hyperbolic orbit such that 
a stable and an unstable curve, say $W^{s}(\overline \varphi_0,\overline\theta_0)$ and
$W^{u}(\overline \varphi_j,\overline\theta_j)$ intersect tangentially at $(\varphi_0,\theta_0)$. 
Then $\alpha$ can be approximated by curves in $U_N$ such that $(\varphi_0,\theta_0)$ is a transverse heteroclinic (or homoclinic) point of
the associated billiard map.
\end{lemma}
\Proof  
Since $(\varphi_0,\theta_0) \in 
\left( W^{s}(\overline \varphi_0,\overline\theta_0) \cap 
W^{u}(\overline \varphi_j,\overline\theta_j) \right)$
there are sequences 
$(\varphi_{i},\theta_{i})=T_\alpha^{in}(\varphi_0,\theta_0) \to (\overline \varphi_0,\overline\theta_0)$ and
$(\varphi_{-i},\theta_{-i})=T_\alpha^{-in}(\varphi_0,\theta_0)\to 
(\overline \varphi_j,\overline\theta_j)$ as 
$i\to+\infty$. Then, there exists an interval $I_1$ such that $\varphi_0\in I_1$, 
$\overline \varphi_0, \overline \varphi_1,...,\overline \varphi_n\notin I_1$ and 
$\varphi_m\notin I_1, \forall m\in Z$.

As $T_\alpha$ is a $C^{l-1}$-diffeomorphism, each $(\varphi_{i},\theta_{i})$ is also a heteroclinic tangent point for every $i$. Moreover, as $T_\alpha$ has the Twist property, it is not possible to the tangency at every $(\varphi_{in},\theta_{in})$ to be vertical. So we can suppose that the stable and unstable
curves are local graphs over the $\varphi$-axis at a neighbourhood of, for instance,
$(\varphi_0,\theta_0)$. So there is an interval $I_2$ containing $\varphi_0$ such that
$W^{s}(\overline \varphi_0,\overline\theta_0)$ and 
$W^{u}(\overline \varphi_j,\overline\theta_j)$ are given locally by the graphs of 
$\theta=\theta^{s}(\varphi)$ and $\theta=\theta^{u}(\varphi)$,
with $\theta^{s}(\varphi_0)=\theta^u(\varphi_0)=\theta_0$ and 
$\frac{d\theta^{s}}{d\varphi}(\varphi_0)=\frac{d\theta^{u}}{d\varphi}(\varphi_0)=\theta'_0$.

Those graphs define two pencils of rays that focuses forward and backward at 
the distances (see, for instance \cite{cmp})
$$d_+=\frac{R_0\sin\theta_0}{1+\theta'_0}\hskip1cm\mbox{and}\hskip1cm
d_-=\frac{R_0\sin\theta_0}{1-\theta'_0},$$
where $R_0=R(\varphi_0)$ is the radius of curvature of $\alpha$ at $\varphi_0$.

Let $I=I_1\cap I_2$ and $\tilde\alpha$ be a normal perturbation 
$\tilde\alpha (\varphi)=\alpha (\varphi)+\lambda(\varphi)\eta(\varphi)$
where $\lambda(\varphi)=0$ if $\varphi\notin I$, $\lambda(\varphi_0)=\lambda'(\varphi_0)=0$, 
$\lambda''(\varphi_0)\neq 0$ and $||\lambda||_2$ is sufficiently 
small in order that $\tilde \alpha\in U_N$ and $T_\alpha$ and $T_{\tilde\alpha}$ are $C^1$-close. The two curves $\alpha$ and $\tilde\alpha$ have a contact of order 1 at $\varphi_0$ and the radius of curvature of $\tilde\alpha$ at this contact point is $\tilde R_0=R_0-\lambda''(\varphi_0)$.
As $\tilde \alpha$ and $\alpha$ differs only on $ I\setminus \{\varphi_0\}$, 
$\{(\overline \varphi_0,\overline \theta_0),...,(\overline \varphi_{n-1},
\overline \theta_{n-1})\}$ is also a $n$-periodic hyperbolic orbit for $T_{\tilde\alpha}$. Moreover, every 
trajectory not hitting $\tilde \alpha (I)$ is the same both for $\alpha$ and for 
$\tilde\alpha$.

Let $\tilde W^{s}(\overline \varphi_0,\overline\theta_0)$ and 
$\tilde W^{u}(\overline \varphi_j,\overline\theta_j)$ be the stable and unstable curves 
of $(\overline \varphi_0,\overline\theta_0)$ and $(\overline \varphi_j,\overline\theta_j)$ 
under $T_{\tilde\alpha}$. We can choose $||\lambda||_2$ sufficiently small such that both are also local graphs over $I$ given by $\theta=\tilde\theta^{u}(\psi)$ and 
$\theta=\tilde\theta^{s}(\psi)$, where $\psi=\psi (\varphi)$ is the angular
parameter of $\tilde \alpha$ with $\varphi_0=\psi(\varphi_0)$.
The pencil of rays defined by $(\psi,\tilde\theta^{u}(\psi))$ will focuses
backward at the distance $\tilde d_-^u$ and the pencil
 $(\psi,\tilde\theta^{s}(\psi))$ will focuses forward at the distance $\tilde d_+^{s}$ with 
$$\tilde d_-^u=\frac{( R_0- \lambda''(\varphi_0))\sin\theta_0}
{1-\frac{d\tilde\theta^u}{d\psi}(\varphi_0)}\hskip1cm\mbox{and}\hskip1cm
\tilde d_+^{s}=\frac{( R_0- \lambda''(\varphi_0))\sin\theta_0}
{1+\frac{d\tilde\theta^{s}}{d\psi}(\varphi_0)}.$$
As the contact of $\alpha$ and $\tilde\alpha$ on $\varphi_0$ is of order 1, preserving the point and the tangent, the trajectory of $(\varphi_0,\theta_0)$ is also the same for both billiards. Then 
$(\varphi_0,\theta_0)$ is a heteroclinic point for $T_{\tilde\alpha}$ and 
$\theta_0=\tilde\theta^{u}(\varphi_0)=\tilde\theta^{s}(\varphi_0)$.

As the curve is unchanged outside $I$, the beam of trajectories given by 
$(\psi,\tilde\theta^{u}(\psi))$ remains the same until it hits $\tilde\alpha (I)$,
implying that $\tilde d_-^u=d_-$ and then 
$$\frac{d\tilde\theta^u}{d\psi}(\varphi_0)=\theta'_0+ \lambda''(\varphi_0)\frac{1-\theta'_0}{R_0}$$

Applying the same construction for the stable curves , with $T_{\alpha}^{-1}$ and $T_{\tilde\alpha}^{-1}$ gives $\tilde d_+^{s}=d_+$
and 
$$\frac{d\tilde\theta^{s}}{d\psi}(\varphi_0)=\theta'_0- \lambda''(\varphi_0)\frac{1+\theta'_0}{R_0}
\neq \frac{d\tilde\theta^u}{d\psi}(\psi_0)$$
implies that the invariant curves for $T_{\tilde\alpha}$ will intersect transversally.\qed

The same reasoning also works for the invariant curves of points on different hyperbolic orbits. As we have only a countable number of hyperbolic orbits, each one with a finite number of points, we can conclude that:
\begin{theorem}
Generically, for billiards on ovals , the invariant curves of two hyperbolic points are transversal.
\end{theorem}

We remark that we do not prove that every homo/heteroclinic orbit is transversal. We do know that generically two invariant stable and unstable curves either do not intersect or have at least one transversal homoclinic orbit, but there can also be tangent orbits.

\section{Generic Dynamics}

\subsection{Rotational Invariant Curves}
A closed, simple, continuous curve $\gamma \subset [0,2\pi)\times(0,\pi)$ which is not homotopic to a point is called a rotational curve.
It is invariant if $T_\alpha(\gamma)=\gamma$. 
The phase-space of the circular billiard, for instance, is foliated by rotational invariant curves.
On the other hand, there are billiards on ovals with no rotational invariant curves at all, as showed in  \cite{hal}.
However, for sufficiently differentiable ovals (\cite{dou}, \cite{laz}),  the Twist property implies that $T_\alpha$ has rotational invariant curves in any small neighbourhood of the boundaries $B_0=[0,2\pi)\times\{0\}$ and $B_\pi=[0,2\pi)\times\{\pi\}$ of the cylinder $[0,2\pi)\times(0,\pi)$.

\begin{prop}
For generic oval billiards the rotation number of any rotational invariant curve is irrational.
\end{prop}
\Proof A rotational invariant curve $\gamma$ is a Lipschitz graph over $[0,2\pi)$ \cite{bir32}. So $\gamma(\varphi)=(\varphi,\theta(\varphi))$, where $\theta$ is continuous and there exists $L>0$ such that $|\theta(\varphi_1)-\theta(\varphi_0)|\leq L|\varphi_1-\varphi_0|$, $\forall \varphi_1, \varphi_0 \in [0,2\pi)$.
Let $\Pi(\varphi,\theta)=\varphi$ and $f(\varphi)=\Pi(T_\alpha(\varphi,\theta(\varphi))$. Then $f$ is a homeomorphism of the circle and so its degree is $\pm 1$.

Suppose that the rotation number of $\gamma$ is rational. If deg$(f)=1$ then $f$ has periodic orbits, all with the same period. If deg$(f)=-1$ then $f^2$ has periodic orbits, all with the same period. 
As $T_\alpha$ is generic, there is only a finite number of periodic orbits, for each period, and they are nondegenerate. Hence there exists an $N$ such that $F=f^N$ is a homeomorphism with a finite number of fixed points, all nondegenerate. Let $\psi_1<\psi_2$ be two consecutive fixed points. Then for every $\varphi\in (\psi_1,\psi_2)$, $\varphi<F(\varphi)<...<F^j(\varphi)<F^{j+1}(\varphi)<...<\psi_2$, $F^j(\varphi)\to \psi_2$ as $j\to\infty$ and $F^j(\varphi)\to \psi_1$ as $j\to-\infty$.
Let $\epsilon>0$ be such that $\psi_1<\psi_2-\epsilon$ and let
$\gamma^-_\epsilon=\{(\varphi,\theta(\varphi)), \psi_2-\epsilon<\varphi\leq\psi_2\}$ and $C^-_\epsilon =\{(\varphi,\theta), \psi_2-\epsilon<\varphi\leq\psi_2, |\theta-\theta(\psi_2)|\leq L|\varphi-\psi_2|
\}$. Clearly $T_\alpha^{jN}(\gamma^-_\epsilon)\subset\gamma^-_\epsilon\subset C^-_\epsilon$.
Let us suppose that $(\psi_2, \theta(\psi_2))$ is an elliptic periodic point of $T_\alpha$. Then $DT_\alpha^N(\psi_2, \theta(\psi_2))$ is a rotation of angle $\beta\neq 0,\pi, 2\pi$. 
There exists $M$ such that $DT_\alpha^{MN}(\psi_2, \theta(\psi_2))(C^-_\epsilon)\cap C^-_\epsilon=\{(\psi_2, \theta(\psi_2))\}$.
As $T_\alpha$ is at least $C^1$ then there is $\epsilon_0\leq \epsilon$ such that
$(\psi_2, \theta(\psi_2))$ is the unique fixed point of $T_\alpha^{MN}$ in $C^-_{\epsilon_0}$ and $T_\alpha^{MN}(C^-_{\epsilon_0})\cap C^-_{\epsilon_0}=\{(\psi_2, \theta(\psi_2))\}$. But this is impossible since 
$T_\alpha^{MN}(\gamma^-_{\epsilon_0})\subset C^-_{\epsilon_0}$.

Then all periodic points in $\gamma$ must be hyperbolic. But this is also impossible in the generic case since $\gamma$ will be the union of the periodic points and saddle connections.
So, the rotation number must be irrational. \qed

From the proof of the above proposition we have the following results:
\begin{lemma}\label{lem:irrac}
Generically, rotational invariant curves can not cross the middle segment $[0,2\pi)\times\{\pi/2\}$.
\end{lemma}
\Proof
Let, as before, $H(\varphi,\theta)=(\varphi,\pi-\theta)$ be the reversing symmetry. Clearly $\gamma$ is a rotational invariant curve if and only if $H(\gamma)$ is also a rotational invariant curve. Suppose that a point  $(\varphi_0,\pi/2)\in \gamma$. Then $(\varphi_0,\pi/2)\in \gamma\cap H(\gamma)$ and, since they are invariant,  $T_\alpha (\varphi_0,\pi/2)=(\varphi_1,\pi/2)\in \gamma\cap H(\gamma)$ and   $\{(\varphi_0,\pi/2),(\varphi_1,\pi/2)\}$ is a 2-periodic orbit in $\gamma$, which is impossible.
Then, for a generic oval billiard, if $(\varphi_0,\theta_0)\in \gamma$ and $\theta_0<\pi/2$ (resp. $\theta_0>\pi/2$) then $\theta(\varphi)<\pi/2$ (resp. $\theta(\varphi)>\pi/2$) for all $\varphi\in [0,2\pi)$ or, in other words, the orbits on $\gamma$ respect the order of the circle $[0,2\pi)$ (respec. reverse the order). \qed 

\begin{lemma}\label{lem:rot}
Let $\gamma$ be a rotational invariant curve and ${\cal O}(p)$ be a hyperbolic periodic orbit, with unstable manifold $W^u$ and stable manifold $W^s$. Then $W^u\cap\gamma=\emptyset$ and $W^s\cap\gamma=\emptyset$.
\end{lemma}
\Proof
Suppose that $x\in W^s\cap\gamma$. Then $T_\alpha^jx\to {\cal O}(p)$ and then ${\cal O}(p)\subset\gamma$, since $\gamma$ is a continuous invariant curve.
But $\gamma$ has irrational rotation number and then do not contain any periodic orbit.
The argument for $W^u$ is analogous. \qed

\subsection{Instability regions}

Birkhoff \cite{bir32} called the region between two invariant rotational curves, with no other invariant rotational curves inside, an instability region. In this subsection we will characterize the instability regions for a generic oval billiard.

\begin{definition}
A cylinder $R\subset [0,2\pi)\times [0,\pi]$ is a non-empty closed connected set such that $\Pi(R)=[0,2\pi)$ and whose boundaries are two continuous rotational curves $r_1$ and $r_2$.
\end{definition}
Note that we are not asking neither the cylinder nor the boundary curves to be invariant under $T_{\alpha}$. 

\begin{theorem}\label{teo:inst}
Let $T_\alpha$ be generic and ${\cal O}(p)$ be a hyperbolic periodic orbit such that its unstable manifold $W^u$ satisfies $\Pi(W^u)=[0,2\pi)$. Then the smallest cylinder containing $\overline {W^u}$ is an instability region.
\end{theorem}

Before proving this theorem, let us remark that for a generic oval billiard and for each period $n$, there is only a finite number of $n$-periodic hyperbolic orbits and at least one of them is a global minimizer of the action $G_{m,n}$, for a suitable $m$. This implies that it has heteroclinic connections \cite{ban} and then $\Pi(W^u)=[0,2\pi)$. Hence the projection hypothesis applies to at least one orbit of each period.

\Proof 
Let $p=(\overline{\varphi_0},\overline{\theta_0})$ be a hyperbolic periodic point, ${\cal O}(p)$ its orbit and $W^u$ the unstable manifold of ${\cal O}(p)$, with $\Pi(W^u)=[0,2\pi)$. 

First, assume that there is an invariant rotational curve $\gamma(\varphi)=(\varphi,\theta(\varphi))$ such that ${\cal O}(p)$ lies in the cylinder $C_\gamma$ bounded by $\gamma$ and $H(\gamma)$.
Since, by lemma~\ref{lem:irrac}, $\gamma$ has irrational rotation number we have that $p\notin\gamma$ and $W^u\cap \gamma=W^u\cap H(\gamma)=\emptyset$. Moreover, $C_\gamma$ is invariant and contains $\overline {W^u}$.

Let then, $C \subset C_{\gamma}$ be the smallest cylinder containing $\overline {W^u}$. Since 
$T_{\alpha} C$ is also a cylinder and $W^u \subset C$,  $T_{\alpha} W^u = W^u\subset T_{\alpha} C$, and we have that $\overline {W^u}\subset T_{\alpha}C$.
So, $T_{\alpha}C$ is also a cylinder containing $\overline {W^u}$.
and as $C$ is the smallest cylinder, we must have $C \subset  T_{\alpha}C$.
However, $T_{\alpha}$ is area-preserving, so we conclude that  $T_{\alpha}C = C$, i.e, $C$ is an  invariant cylinder.

As $T_{\alpha}$ is a diffeomorphism,  $T_{\alpha}\partial C=\partial C$ and $T_{\alpha} int(C)=int(C)$. From Birkhoff's Theorem it follows that $\partial C$ is the union of two rotational invariant curves $c_1$ and $c_2$.

In fact, these two curves $c_1$ and $c_2$ are contained in $\overline {W^u}$. To prove this fact, let us suppose that there is a point $p\in c_1, p\notin \overline {W^u}$.
Then there is an open ball $B$, centered at $p$ such that $\overline B\cap  \overline {W^u}=\emptyset$. We can then construct a new cylinder ${\tilde C} \subset C$ with $\partial {\tilde C}=\{\tilde c_1,c_2\}$ where $\tilde c_1=c_1$ outside $\overline B$ and $\tilde c_1= \partial \overline B\cap C$ otherwise.
\begin{center}
\includegraphics[width=5cm]{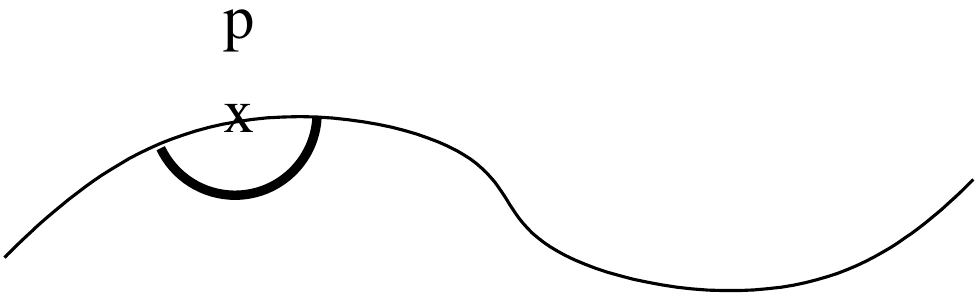}
\end{center}
We have then a new cylinder satisfying $\overline {W^u} \subset {\tilde C}\subset C$ and so $C$ is not the smallest one, which is a contradiction.
Now, as $W^u$ is the union of continuous curves, its closure contains the boundaries $c_1$ and $c_2$ and $W^u\cap\gamma=\emptyset$ for every $\gamma$ rotational and invariant, it follows that there is no invariant curve in $int(C)$.

We have proved that $C$ is an invariant cylinder, its boundaries are rotational invariant curves and there are no other rotational invariant curves inside it. So, it is the instability region containing $\overline {W^u}$.

Suppose now that there is no invariant rotational curves curves bounding $W^u$ from bellow, i.e,
$0<\overline{\theta_0}<\theta(\varphi_0)<\pi-\theta(\varphi_0)$ for all rotational invariant curves $\gamma(\varphi)=(\varphi,\theta(\varphi))$.
We can extend $T_\alpha$ continuously to the circle $B_0=[0,2\pi)\times\{0\}$, observing that $T_\alpha(B_0)=B_0$. Now, any cylinder with boundaries $B_0$ and a rotational invariant curve contains $\overline {W^u}$.  The smallest one, $C$, will have boundaries $B_0$ and $c_2$, because if there were a new boundary $c_1$ then it will be a rotational invariant curve with $\theta(\varphi_0)<\overline{\theta_0}<\pi-\theta(\varphi_0)$.
Reasoning as above we conclude that $C$ is the instability region that contains $\overline {W^u}$.

The case $\theta(\varphi_0)<\pi-\theta(\varphi_0)<\overline{\theta_0}$ for all rotational invariant curves follows from the reversibility property of $T_\alpha$.

Finally, if $T_{\alpha}$ has no rotational invariant curves, the instability region containing $\overline {W^u}$ is the entire cylinder $[0,2\pi)\times [0,\pi]$.
\qed

Obviously, we can make a similar construction and obtain  the same result for the stable manifold $W^s$.

\subsection{Islands of Stability}

\begin{prop} \label{prop:ilha}
Let $V$ be the complement of $\overline {W^u}$ in the instability region $C$.
Then ${V=\cup_1^\infty U_i}$, where all the $U_i$'s are disjoint open sets, homeomorphic to discs. Dynamically, each $U_i$ returns over itself.
\end{prop}
\Proof 
$V$ is open and is the union of its disjoint open connected components. Let $U$ be one of them. $U$ can not divide $C$ into two disjoint sets because $\overline{W^u}$ contains $W^u$, which projects over $[0,2\pi)$ and converges to the two rotational curves of $\partial C$.
So $U$ is homotopically trivial. 

As $V\subset C\subset [0,2\pi)\times [0,\pi]$ has finite area
and each connected  component has positive area. 
So there can not exist a non countable number of them. 

Clearly, $T_{\alpha}V=V$. As $T_{\alpha}$ is area preserving then, given $U_i$, there exists a smallest $n(i)$ such that $T_{\alpha}^{n(i)}U_i\cap U_i\neq \emptyset$. But, as each $U_i$ is a connected component $V$, we must have $T_{\alpha}^{n(i)}U_i= U_i$. In other words, each component $U_i$ of $V=C - \overline{W^{u}}$ is periodic. \qed

The above proposition suggests the definition
\begin{definition}
Given $i$, let $n(i)$ be the smallest integer such that $T_{\alpha}^{n(i)}U_i= U_i$. We call $U_i$ an island of stability and $n(i)$ its period. 
The invariant set $I=U_i\cup T_{\alpha}U_i\cup...\cup T_{\alpha}^{n(i)-1}U_i$ will be called an archipel.
\end{definition}

Remark that, as each island $U_i$ is invariant by $T_\alpha^{n(i)}$ and homeomorphic to a ball then
\begin{prop}
The boundary $\partial U_i$ is a closed connected set, with empty interior, contained in $\overline {W^u}$ and invariant by $T_\alpha^{n(i)}$.
\end{prop}

An easy consequence of proposition \ref{prop:ilha} is
\begin{prop} Every periodic point in an island of stability has a period multiple of the period of the island.
\end{prop}

More strongly, we have that
\begin{prop} Each island of stability contains a periodic point with the same period as the island.
\end{prop}
\Proof
Let $U_i$ be an island of stability with period $n(i)$.
Since the restriction $T^{n(i)}_\alpha\vert U_i$ is area preserving and conjugated to a diffeomorphism of the plane which is orientation preserving,
it follows from Brouwer's Translation Theorem that  $T^{n(i)}_\alpha\vert U_i$ has a fixed point.
So there exists a point $p_i\in U_i$ such that $T_{\alpha}^{n(i)}p_i=p_i$ and so a periodic orbit ${\cal O} (p_i)=\{p_i, T_{\alpha}p_i,..., T_{\alpha}^{n(i)-1}p_i\}$ with the smallest possible period. \qed

Remark that, as billiards have no fixed points then $T_\alpha U_i\cap U_i=\emptyset$ and $n(i)\geq 2$.

\begin{prop} All points in an island of stability have the same rational rotation type.
\label{prop:rotnum}\end{prop}
\Proof
Let, as above, $U_i$ be an island of stability with period $n(i)$. Following \cite{lec87} we observe that if $\hat{U_i} $ is the lift of $U_i$  to the universal cover and $\hat{T_\alpha}$ is a lift of $T_\alpha$, then there exists an integer $m$ such that ${\hat{T_\alpha}}^{n(i)}(\hat{U_i})= \hat{U_i} + (2m\pi,0)$.
Since $\hat{U_i} $ is bounded, this implies that that for all $z \in U_i$ if 
$\hat{z}$ denotes a lift of $z$, then there exists $\displaystyle \lim_{n\rightarrow \infty} \frac{pr_{1}\circ{\hat{T_\alpha}}^{kn(i)}(\hat z )}{k}=\frac{2m\pi}{n(i)}$, where $pr_{1}$ is the projetion on the first factor. \qed

\subsection{The Instability Set}

For the sake of completeness, we conclude by describing some well known dynamical consequences of the area preserving and the twist properties.

Let $T_\alpha$ be generic, in the sense of the previous sections, and ${\cal O}(p)$ be a hyperbolic periodic orbit such that its unstable manifold $W^u$ satisfies $\Pi(W^u)=[0,2\pi)$. Let $C$ be the instability region containing $\overline {W^u}$.
Following the terminology of Franks in \cite{fra} we call $\overline{W^{u}}$ the instability set. 

Clearly, the instability set $\overline{W^u}$ does not contain any Moser stable periodic orbit, since they have open islands around them. For hyperbolic periodic orbits we have

\begin{prop}
Let ${\cal O}(q)$ be a hyperbolic periodic orbit in $\overline{W^u}$ and let $W^s(q)$ and $W^u(q)$ be its invariant stable and unstable manifolds. 
Then either $\overline{W^u(q)}$ and $\overline{W^s(q)}\subset \overline{W^u}$ or $q\in\partial U_i$, ie, the boundary of an island.
Moreover, for $j=u$ or $s$, if $\overline{W^j(q)}\subset \overline{W^u}$ and $\Pi(W^j(q))=[0,2\pi)$ then $\overline{W^j(q)}= \overline{W^u}$;
or if $q\in\partial U_i$ then two branches of the invariant curves of $W^{u}(q)$ and/or $W^{s}(q)$ are contained in $\partial U_i$.
\end{prop}
\Proof
Suppose that there exists $z\neq q$ such that $z\in W^j(q)\cap\overline{W^{u}}$, $j=u,s$. Then by lemma 3.1 of \cite{lec87}, the branch of $W^j(q)$ which contains $z$ is entirely contained in $\overline{W^{u}}$. So if each branch of $W^j(q)$ intersects $\overline{W^{u}}$ then  $W^j(q)\subset\overline{W^{u}}$ and obviously, $\overline{W^j(q)}\subset\overline{W^{u}}$.

Let us now suppose that a branch $\gamma$ of $W^j(q)$ is not in $\overline{W^{u}}$. Then $\gamma$ is contained in an island of stability $U_i$ or is outside the region of instability $C$. As the $\omega,\alpha$-limit of $\gamma$ is $q$, either $q\in\partial U_i$ or $q\in c_1$ or $c_2$, the boundary curves of $C$. But this last case is impossible in the generic case. So $q\in\partial U_i$.
But in this case, it is enough to remark that $\partial U_i$ is an invariant closed connected set and that the only sets with these properties containing $q$ are the invariant manifolds to conclude that two branches of the invariant curves of $W^{u}(q)$ and/or $W^{s}(q)$ are contained in $\partial U_i$.

If $\Pi(W^j(q))=[0,2\pi)$ then $q\notin \partial U_i$ and, by theorem \ref{teo:inst},
$\overline{W^j(q)}= \overline{W^u}$.  \qed

Concerning non-periodic orbits and as a consequence of proposition \ref{prop:rotnum} we have that
\begin{prop}
Any  Aubry- Mather set of irrational type contained in $C$ is actually contained in the instability set $\overline{W^{u}}$.
\end{prop}

And as it was proven by Le Calvez \cite{lec87}
The instability set $\overline{W^{u}}$ contains the closure of the orbits which $\alpha$ and/or $\omega$-limits are $c_1$ and $c_2$.

\section{Final remarks}

So far we have obtained a rough description of the dynamics in the instability region, very similar to what hapens for generic area preserving twist maps.
It is made of the closure of an hyperbolic orbit and the union of periodic islands.

In order to proceed with this description it is necessary to address two basic related themes.
The first one has to do with the islands: the existence of a finite number of them, the existence of a lower bound for the period of an island (e.g. there are examples with no period two islands) and the existence of connecting orbits between the boundaries of the islands. 
A very interesting question is if it is possible to have an instability region without any islands.

The other concerns the properties of $\overline{W^{u}}$: does it have empty interior? Is it topologically transitive?

A starting point is the instability region that contains the period two orbits. In this region, the reversibility symmetry of the phase space, together with geometric properties of the boundary of the billiard may allow us to obtain new results about the above questions.

\bigskip

{\bf Acknowledgments}:
The authors thank Conselho Nacional de Desenvolvimento
Cient\'{\i}fico (CNPq) and Funda\c c\~ao de Amparo a Pesquisa de
Minas Gerais, Brazilian agencies, for financial support. SPC
gratefully acknowledges the hospitality of Laboratoire Emile Picard, Universit\'e Paul Sabatier (Toulouse III), where part of this work was done, under the financial support of Coordena\c c\~ao de Aperfei\c coamento de Pessoal de N\'\i vel Superios (CAPES).


\begin {thebibliography} {99}

\bibitem{ban} V.Bangert: Mather sets for twist maps and geodesics on tori, Dynamics reported vol 1, 1-57, 1988.
\bibitem{bir27}  G.D.Birkhoff: {\em Dynamical Systems}.  
 Providence, RI: A. M. S. Colloquium Publications, 1966, (Original ed. 1927,)
\bibitem{bir32} G.D.Birkhoff: Sur quelques courbes ferm\'ees remarquables, Bull. SMF, {\bf 60}, 1-26, 1932. (Also in {\em Collected Math. Papers of G; D; Birkhoff}, vol II, 444-461.)
\bibitem{etds} M.J.Dias Carneiro, S.Oliffson Kamphorst,S.Pinto de Carvalho:
Elliptic Islands in Strictly Convex Billiards, Erg.Th.Dyn.Sys., {\bf 23/3}, 799-812, 2003.
\bibitem{don} V.J.Donnay: Creating transverse homoclinic connections in planar billiards, Jr. Math. Sci., {\bf 128/2},  2747-2753, 2005.
\bibitem{dou} 
R.Douady: Applications du th\'eor\`eme des tores invariants. Th\`ese de 3\`eme Cycle, Univ. Paris VII, 1982.
\bibitem{fra} J.Franks: Rotation Numbers and Instability Sets, Bull. AMS New Series), {\bf 40/3}, 263-279, 2003.
\bibitem{hal} B.Halpern: Strange Billiard Tables, Trans. Amer. Math. Soc. {\bf 232}, 297-305, 1977.
\bibitem{kat} A.Katok, B.Hasselblat: {\em Introduction to the Modern Theory of Dynamical Systems}. Univ. Cambridge Press, 1997.
\bibitem{koz} V.V.Kozlov, D.V.Treshchev: {\em Billiards - A Genetic Introduction to the Dynamics of Systems with Impacts}, Transl.Math.Monog., AMS, vol 89, 1980.
\bibitem{kozl}
V.V.Kozlov: Two-link billiard trajectories: extremal properties and stability, J.Appl. Maths Mechs, {\bf 64/6}, 903-907, 2000.
\bibitem{laz}
V.F.Lazutkin: The existence of caustics for a billiard problem in a convex domain, Math.USSR Izvestija, {\bf 7/1}, 185-214, 1973.
\bibitem{lec87}
P.Le Calvez: Propri\'et\'es dynamiques des r\'egions d'instabilit\'e, Ann.Sci. Ec. Norm. Sup., tomme 20, 443-464, 1987.
\bibitem{lev}
P.Levallois, M.B.Tabanov: S\'eparation des s\'eparatrices du billard elliptique
pour une p\'erturbation alg\'ebrique et symm\'etrique de l'\'ellipse, C.R.Ac.Sci.Paris
-Sr I-Math, {\bf 316/3}, 589-592, 1993.
\bibitem{mac} R. MacKay, J. Meiss: Linear stability of periodic orbits in lagrangian systems, Phys.Lett. A, {\bf 98}, 92-94, 1983.
\bibitem{cmp} R.Markarian, S.Oliffson Kamphorst e S.Pinto de Carvalho: 
Chaotic properties of the elliptical stadium, Comm. Math. Phys. {\bf 174}, 661-679, 1996.
\bibitem{rag}
C.G.Ragazzo, M.J.Dias Carneiro, S.Addas Zanata: Introdu\c c\~ao \`a Din\^amica de Aplica\c c\~oes do Tipo Twist. Pub.Mat., IMPA, 2005.
\bibitem{tab} S.Tabachnikov: {\em Billiards}, S M F Panorama et Syntheses 1, 1995.

\end {thebibliography}


\end{document}